\documentclass{article}[12pt]

\usepackage{latexsym}
\usepackage{amsmath}
\usepackage{amssymb}

\begin{document}

\title{A class of solutions of the asymmetric May-Leonard model}

\author{Francesco Calogero$^{a,b}$\thanks{e-mail: francesco.calogero@roma1.infn.it}
\thanks{e-mail: francesco.calogero@uniroma1.it}
 , Farrin Payandeh$^c$\thanks{e-mail: farrinpayandeh@yahoo.com}
 \thanks{e-mail: f$\_$payandeh@pnu.ac.ir}}

\maketitle   \centerline{\it $^{a}$Physics Department, University of
Rome "La Sapienza", Rome, Italy}

\maketitle   \centerline{\it $^{b}$INFN, Sezione di Roma 1}

\maketitle

\maketitle   \centerline{\it $^{c}$Department of Physics, Payame
Noor University, PO BOX 19395-3697 Tehran, Iran}

\maketitle

\begin{abstract}

The asymmetric May-Leonard model is a prototypical system of $3$ nonlinearly
coupled first-order Ordinary Differential Equations with second-degree
polynomial right-hand sides. In this short paper we identify a class of
\textit{special} solutions of this system which do not seem to have been
previously advertised in spite of their rather elementary character.

\end{abstract}

\section{Introduction and results}

The May-Leonard model \cite{ML1975} is a prototypical system, introduced
almost half a century ago, of $3$ nonlinearly coupled first-order Ordinary
Differential Equations with second-degree polynomial right-hand sides. In
self-evident notation it reads as follows:
\begin{subequations}
\label{ML1975}
\begin{equation}
\dot{\xi}_{1}\left( t\right) =\xi _{1}\left( t\right) \left[ 1-\xi
_{1}\left( t\right) -\alpha \xi _{2}\left( t\right) -\beta \xi _{3}\left(
t\right) \right] ~,
\end{equation}%
\begin{equation}
\dot{\xi}_{2}\left( t\right) =\xi _{2}\left( t\right) \left[ 1-\beta
x_{1}\left( t\right) -x_{2}\left( t\right) -\alpha x_{3}\left( t\right) %
\right] ~,
\end{equation}%
\begin{equation}
\dot{\xi}_{3}\left( t\right) =\xi _{3}\left( t\right) \left[ 1-\alpha
x_{1}\left( t\right) -\beta x_{2}\left( t\right) -x_{3}\left( t\right) %
\right] ~.
\end{equation}

The \textit{asymmetric} May-Leonard model \textbf{(}much investigated in the
literature\textbf{: }see for instance \cite{CSL1998} \cite{AFRS2019} \cite%
{LMV2020}, and references therein) is the following generalization---again,
in self-evident notation---of the May-Leonard model (\ref{ML1975}):
\end{subequations}
\begin{subequations}
\label{ML}
\begin{equation}
\dot{x}_{1}\left( t\right) =x_{1}\left( t\right) \left[ \eta -x_{1}\left(
t\right) -a_{12}x_{2}\left( t\right) -a_{13}x_{3}\left( t\right) \right] ~,
\label{ML1}
\end{equation}%
\begin{equation}
\dot{x}_{2}\left( t\right) =x_{2}\left( t\right) \left[ \eta
-a_{21}x_{1}\left( t\right) -x_{2}\left( t\right) -a_{23}x_{3}\left(
t\right) \right] ~,  \label{ML2}
\end{equation}%
\begin{equation}
\dot{x}_{3}\left( t\right) =x_{3}\left( t\right) \left[ \eta
-a_{31}x_{1}\left( t\right) -a_{32}x_{2}\left( t\right) -x_{3}\left(
t\right) \right] ~.  \label{ML3}
\end{equation}%
Here and hereafter the $6$ (\textit{constant}) parameters $a_{nm}$ ($n=1,2;$
$m=1,2,3$) are \textit{a priori arbitrary} (except, of course, for the
restrictions on their values identified below); and superimposed dots
indicate differentiations with respect to the independent variable $t$
("time").

\textbf{Remark 1}. The original May-Leonard model \textbf{(}\ref{ML1975})
clearly corresponds---up to a trivial renaming of the dependent
variables---to the asymmetric May-Leonard model (\ref{ML}) with
\end{subequations}
\begin{equation}
a_{12}=a_{23}=a_{31}=\alpha ~;~~~a_{13}=a_{21}=a_{32}=\beta ~,~~~\eta
=1~.~~~\blacksquare
\end{equation}

\textbf{Remark 2}. Often the asymmetric May-Leonard model is characterized
by the system of $3$ ODEs (\ref{ML}) with $\eta =1$. It is of course trivial
to reduce the system (\ref{ML}) to its more standard version with $\eta =1$
via the following simultaneous rescaling of the independent variable $t$ and
the $3$ dependent variables $x_{n}\left( t\right) $: $x_{n}\left( t\right)
\Rightarrow \eta \tilde{x}_{n}\left( \tilde{t}\right) ,$ $\tilde{t}=\eta t$.
We prefer to keep the extra parameter $\eta $ in view of its relevance to
generate the \textit{isochronous} variant of the asymmetric May-Leonard
model, see below \textbf{Remark 4}. \ \ $\blacksquare $

Via the following well-known (see, for instance, \cite{FC2008}) change of
variables,
\begin{subequations}
\begin{equation}
x_{n}\left( t\right) =\exp \left( \eta t\right) y_{n}\left( \tau \right)
~,~~~\tau =\left[ \exp \left( \eta t\right) -1\right] /\eta ~,~~~n=1,2,3
\label{xnyn}
\end{equation}%
implying%
\begin{equation}
x_{n}\left( 0\right) =y_{n}\left( 0\right) ~,  \label{xnyn0}
\end{equation}%
the new dependent variables $y_{n}\left( \tau \right) $ are easily seen to
satisfy the following (still \textit{autonomous}) system of $3$ nonlinearly
coupled ODEs:
\end{subequations}
\begin{subequations}
\label{ynprime}
\begin{equation}
y_{1}^{\prime }\left( \tau \right) =-y_{1}\left( \tau \right) \left[
y_{1}\left( \tau \right) +a_{12}y_{2}\left( \tau \right) +a_{13}y_{3}\left(
\tau \right) \right] ~,  \label{y1prime}
\end{equation}%
\begin{equation}
y_{2}^{\prime }\left( \tau \right) =-y_{2}\left( \tau \right) \left[
a_{21}y_{1}\left( \tau \right) +y_{2}\left( \tau \right) +a_{23}y_{3}\left(
\tau \right) \right] ~,  \label{y2prime}
\end{equation}%
\begin{equation}
y_{3}^{\prime }\left( \tau \right) =-y_{3}\left( \tau \right) \left[
a_{31}y_{1}\left( \tau \right) +a_{32}y_{2}\left( \tau \right) +y_{3}\left(
\tau \right) \right] ~.  \label{y3prime}
\end{equation}%
Here and below appended primes indicated of course differentiations with
respect to the independent (generally \textit{complex}) variable $\tau $.

The fact that the $3$ polynomials in the right-hand sides of this \textit{%
autonomous} system of $3$ ODEs are \textit{homogeneous} in the dependent
variables $y_{n}\left( \tau \right) $ implies (see \cite{CP2021}) that this
system features the following elementary solution of its initial-values
problem:
\end{subequations}
\begin{equation}
y_{n}\left( \tau \right) =y_{n}\left( 0\right) /\left( 1+z\tau \right)
~,~~~n=1,2,3~,  \label{yntauSol}
\end{equation}%
provided the parameter $z$, the $3$ initial-values $y_{n}\left( 0\right) $
and the $6$ parameters $a_{nm}$ ($n=1,2;$ $m=1,2,3$) satisfy the following $%
3 $ simple algebraic relations:
\begin{subequations}
\label{zz}
\begin{equation}
z=\left[ y_{1}\left( 0\right) +a_{12}y_{2}\left( 0\right) +a_{13}y_{3}\left(
0\right) \right] ~,  \label{z1}
\end{equation}%
\begin{equation}
z=\left[ a_{21}y_{1}\left( 0\right) +y_{2}\left( 0\right) +a_{23}y_{3}\left(
0\right) \right] ~,  \label{z2}
\end{equation}%
\begin{equation}
z=\left[ a_{31}y_{1}\left( 0\right) +a_{32}y_{2}\left( 0\right) +y_{3}\left(
0\right) \right] ~.  \label{z3}
\end{equation}

We can therefore formulate (via (\ref{xnyn})) our main finding.

\textbf{Proposition}: the asymmetric May-Leonard model (\ref{ML}) admits the
following solution of its initial-values problem,
\end{subequations}
\begin{equation}
x_{n}\left( t\right) =x_{n}\left( 0\right) /\left\{ \exp \left( -\eta
t\right) +z\left[ 1-\exp \left( -\eta t\right) \right] /\eta \right\}
~,~~~n=1,2,3~,  \label{xnt}
\end{equation}%
with%
\begin{eqnarray}
z &=&x_{1}\left( 0\right) +a_{12}x_{2}\left( 0\right) +a_{13}x_{3}\left(
0\right)  \nonumber \\
&=&a_{21}x_{1}\left( 0\right) +x_{2}\left( 0\right) +a_{23}x_{3}\left(
0\right)  \nonumber \\
&=&a_{31}x_{1}\left( 0\right) +a_{32}x_{2}\left( 0\right) +x_{3}\left(
0\right) ~,  \label{z}
\end{eqnarray}%
provided the $3$ initial-values $x_{n}\left( 0\right) $ and the $6$
parameters $a_{nm}$ ($n=1,2;$ $m=1,2,3$) satisfy the following $2$ simple
algebraic relations (implied by (\ref{z})):
\begin{eqnarray}
&&x_{1}\left( 0\right) +a_{12}x_{2}\left( 0\right) +a_{13}x_{3}\left(
0\right)  \nonumber \\
&=&a_{21}x_{1}\left( 0\right) +x_{2}\left( 0\right) +a_{23}x_{3}\left(
0\right)  \nonumber \\
&=&a_{31}x_{1}\left( 0\right) .+a_{32}x_{2}\left( 0\right) +x_{3}\left(
0\right) ~.~~~\blacksquare  \label{2Eqs}
\end{eqnarray}

\textbf{Remark 3}. For \textit{any given} assignment of any subset of $7$
out of the $9$ parameters $a_{12},$ $a_{13},$ $a_{21},$ $a_{23},$ $a_{31},$ $%
a_{32},$ $x_{1}\left( 0\right) ,$ $x_{2}\left( 0\right) ,$ $x_{3}\left(
0\right) ,$ these relations (\ref{2Eqs}) determine (easily and uniquely) the
other $2$ of these $9$ parameters in terms of the $7$ \textit{arbitrarily}
assigned; thereby identifying---via (\ref{xnt}) with (\ref{z})---the
corresponding \textit{explicit} solution of the initial-values problem of
the asymmetric May-Leonard model (\ref{ML}).~~~$\blacksquare $

\textbf{Remark 4}. It is obvious that if $\eta $ is an \textit{imaginary}
number, $\eta =\mathbf{i}\omega $ (with $\mathbf{i}$ the \textit{imaginary}
unit, $\mathbf{i}^{2}=-1$, and $\omega $ an \textit{arbitrary nonvanishing
real} number), then \textit{all} the solutions of the asymmetric May-Leonard
model (\ref{ML}) given by (\ref{xnt}) with (\ref{z}) are \textit{isochronous}%
, i. e. \textit{completely periodic }with period $T=2\pi /\left\vert \omega
\right\vert $, $x_{n}\left( t+T\right) =x_{n}\left( t\right) $. However in
this case one would be dealing with a model involving \textit{complex}
dependent variables, $x_{n}\left( t\right) = $Re$ \left[ x_{n}\left(
t\right) \right] +\mathbf{i} $Im$ \left[ x_{n}\left( t\right) \right] ,$
i. e. $6$ \textit{real} variables rather than only $3$; and it would be then
reasonable to also \textit{double} the number of \textit{real} parameters by
setting $a_{nm}= $Re$ \left[ a_{nm}\right] +\mathbf{i} $Im$ \left[
a_{nm}\right] $. And we leave to the interested reader to consider the
behavior of the solution (\ref{xnt}) in the more general case when the
parameter $\eta $ is itself a \textit{complex} number, $\eta = $Re$ \left[
\eta \right] +\mathbf{i} $Im$ \left[ \eta \right] $. \ \ $\blacksquare $

\bigskip

\textbf{Acknowledgements}. It is a pleasure to thank our colleagues Robert
Conte, Fran\c{c}ois Leyvraz and Andrea Giansanti for very useful
discussions. We also like to acknowledge with thanks $2$ grants,
facilitating our collaboration---mainly developed via e-mail exchanges---by
making it possible for FP to visit twice the Department of Physics of the
University of Rome "La Sapienza": one granted by that University, and one
granted jointly by the Istituto Nazionale di Alta Matematica (INdAM) of that
University and by the International Institute of Theoretical Physics (ICTP)
in Trieste in the framework of the ICTP-INdAM "Research in Pairs" Programme.%
\textbf{\ }Finally, we also\ like to thank Fernanda Lupinacci who, in these
difficult times---with extreme efficiency and kindness---facilitated all the
arrangements necessary for the presence of FP with her family in Rome.

\end{document}